\documentclass[letterpaper,10pt,conference]{ieeeconf}
\IEEEoverridecommandlockouts
\usepackage[utf8]{inputenc}
\usepackage{graphicx}
\usepackage[hidelinks]{hyperref}
\usepackage{amsmath}
\allowdisplaybreaks[2]
\usepackage{amsfonts}
\usepackage{soul}
\usepackage{xcolor}
\usepackage{cancel}

\def\R{{\mathbb{R}}}
\def\N{{\mathbb{N}}}

\def\tt{{t}}

\newtheorem{theorem}{Theorem} 
\newtheorem{lemma}{Lemma} 
\newtheorem{fact}{Fact}
\newtheorem{definition}{Definition} 
\newtheorem{remark}{Remark} 
\newtheorem{assumption}{Assumption}
\newtheorem{example0}{Example} 

\title{
Conditions for Normality of
Pontryagin's Maximum Principle 
in the Absence of State Constraints}

\author{
\authorblockN{Andr{\'e} L.~Tits}
\authorblockA{ECE\\
\& Institute for Systems Research\\
University of Maryland\\
College Park, MD 20742, USA\\
(\texttt{andre@umd.edu})\\
ORCID: 0000-0002-1542-7889.\\
Corresponding author.}
\and
\authorblockN{Pierre-Antoine Absil}
\authorblockA{ICTEAM Institute, UCLouvain\\
1348 Louvain-la-Neuve, Belgium\\
ORCID: 0000-0003-2946-4178}
\and
\authorblockN{No\'e Tits}
\authorblockA{ISIA Lab, Numediart Institute\\
UMONS\\
7000 Mons, Belgium\\
ORCID: 0000-0002-1971-4412}
}
\date{\today}

\begin{document}

\maketitle

\begin{abstract}
Pontryagin’s Maximum Principle is said to be ``normal'' for the problem
at hand when it 
admits a maximizing control
that explicitly involves
the problem's cost function, making the maximum principle helpful (while
in the opposite situation it is of no help).
The aim of the present paper is, for Lagrange optimal control problems
without state constraints,
to identify classes of problems for which the Principle is
normal.  General conditions are introduced that guarantee normality,
and classes of problems that satisfy these conditions are identified.
Such classes include problems for which the constrained (by bounds on the control values) ``augmented'' reachable set, with cost along an upward vertical axis, has a convex lower-boundary function and the optimal terminal state lies in the relative interior of the reachable-state set, with the terminal constraint within the affine hull of that set.
Such classes include problems for which the constrained (by bounds on
the values of the control signal) ``augmented'' (with the cost value along
an upward “vertical axis”) reachable set has a “convex” lower boundary and
the optimal terminal state lies in the relative interior of the set
of reachable states.  Special attention is given to Lagrange {\em steering}
problems: there, sufficient conditions are obtained on a nested sequence of
sub-classes of such problems, culminating with the class of steering
problems with scalar control signal and linear time-invariant (LTI) state
dynamics and objective integrand, for which normality is guaranteed whenever all
eigenvalues of the (real) state dynamics matrix $A$ are real (e.g., if $A$ 
is symmetric or triangular or the system is overdamped), but that is shown
to not extend to all 
real $A$ matrices;
still, at least in the $2\times 2$ case, in the presence 
of non-real eigenvalues, normality does holds as long as the imaginary parts
are small enough.

\end{abstract}

\begin{keywords}
optimal control, Pontryagin's maximum principle (PMP), normality of PMP, steering problems
\end{keywords}
\section{Introduction}

Normality of Pontryagin's maximum principle (PMP) has 
been an object of much attention over the past 
few decades; see, e.g.,~\cite{FrankowskaTonon2013, BettiolFrankowska2007,
RampazzoVinter1999, BettiolKhalil2016, FerreiraFontesVinter1999}.
In particular, in the presence of both endpoints constraints, sufficient 
conditions
for normality were obtained in~\cite{FrankowskaTonon2013} that apply
to a large class of optimal control problems with (or without) state 
constraints.  However, as noted just before Theorem~5.6 in that paper,
such conditions are obtained only for trajectories ending on the boundary
$\partial K$ of the set $K$ that defines the state constraints (see (1.2)
in~\cite{FrankowskaTonon2013}); accordingly, in the absence of state
constraints (there are none in the present paper), $K$ is the 
entire space so $\partial K$ is empty, and no trajectories can end 
in $\partial K$, so that Theorem~5.6 
of~\cite{FrankowskaTonon2013} is not relevant when both endpoints
are constrained.

Standard statements of PMP's necessary
conditions of optimality 
(see, e.g.,~\cite{Varaiya72:0, Vinter2000, Liberzon2012}) 
involve a 
constant scalar ``multiplier'' associated to the objective function, 
sometimes
referred to as the abnormal multiplier---noted $p_0^*$ in 
the present paper---for which, depending on the specifics of the problem at 
hand, the value {\em zero}
may be the only possible choice, making the 
PMP abnormal.
The intent of the present paper is to identify classes of optimal
control problems for which PMP's necessary conditions do not 
suffer such drawback.

With the focus first placed on the Mayer problem,
after setting the background
(section~\ref{sec:background}), in section~\ref{sec:p0NOT0},
we introduce conditions (TS1 and TS2) at the terminal state
(LTS1 in a Lagrange-focused context; see~section~\ref{sec:Lagrange Problem} 
below), which are shown
(Theorem~\ref{thm:abnormal}) to be sufficient for PMP's 
conditions to hold with $p_0^*\not=0$ (equivalently, via scaling,
with $p_0^*=1$) 
for the problem under consideration.

The main portion of the paper (section~\ref{sec:Lagrange Problem}) focuses 
on the {\em Lagrange} formulation of optimal control problems.  
It investigates conditions under which TS1/LTS1 and TS2 hold.
(LTS1 excludes downwards directions in the closed convex hull of the 
augmented reachable set's contingent cone.
%

It is then first shown (section~\ref{sec:LTI dynamics}) that, 
for rather general Lagrange problems, 
convexity of the lower boundary function
and an interior point assumption 
imply LTS1, provided the terminal
condition is regular within the affine hull of the reachable set.
This is then shown 
to be the case in a class of problems studied by
L.W.~Neustadt~\cite[Theorem 1]{Neustadt1963}.  

Attention then shifts (section~\ref{sec:LTIstate}) to the class of 
Lagrange {\em steering} problems, specifically,
problems in which the terminal state is prescribed; for such problems,
Condition TS2 is {\em always} satisfied; further, when the state
dynamics is linear and time-invariant (LTI), for a large class of
cost integrands linear in the state, either (i) the optimal control is
bang--bang or (ii) PMP is normal (Theorem~\ref{col:xfin K(T_0)}).
Finally (section~\ref{sec:LTIdynANDcost}), when in addition the cost 
integrand is time-independent and linear 
in state and control signals, PMP is normal
as long as (Theorem~\ref{prop:vertical tangent}) all eigenvalues 
of the state dynamics matrix $A$ are real 
(Theorem~\ref{prop:vertical tangent})---as well as, at least in the 
$2\times 2$ case, when the imaginary parts of the eigenvalues of $A$ 
are appropriately small.

From section~\ref{sec:LTI dynamics} onward, the analysis is limited to
the single-input case; the multi-input case is left for future studies.

\smallskip
{\bf Notation.}  It is mostly standard.  Worth mentioning are the following:
$0_k$ is the origin of $\R^k$; given a scalar $\delta>0$,
$B_\delta(\xi)$ is the open ball of radius $\delta$ centered at $\xi$;
given a differentiable mapping $G:\R^p\to\R^q$, J$G:\R^p\to\R^{q\times p}$
is its Jacobian; and the sign function ``sgn'' takes values 
in $\{-1,0,1\}$, with sgn(0):=0. Further, given $v\in\R^p$ and $w\in\R^q$ the
Matlab notation $[v;w]\in\R^{p+q}$ designates the concatenation of $v$ and $w$,
${\mathcal N}(M)$ and ${\mathcal R}(M)$ denote the kernel (nullspace)
and range of matrix $M$; and,
given a set $S$, $\partial S$ denotes its boundary, co($S$) its convex hull,
cl($S$) is closure, aff($S$) its affine hull, and relint($S$) its relative 
interior.

\section{Background}
\label{sec:background}

Consider the Mayer problem with fixed initial and terminal times,
fixed initial state, and free terminal state, viz.
\begin{equation}
\label{eq:MayerInt}
\begin{aligned}
\min_{u\in\mathcal U}&\ \psi(x(1))\\
\text{s.t.~}&\ \dot x(t)=f(t,x(t),u(t))~\text{a.e.~}t\in[0,1],\\
&\ x\text{~abs.~cont.},~x(0)=\xi^0,
\end{aligned}
\end{equation}
where $\psi:\R^n\to\R$ is continuously differentiable, 
$f:\R\times\R^n\times\R^m\to\R^n$ is (jointly) continuous,
as well as locally Lipschitzian in its second argument (as is
$f_0:\R\times\R^n\times\R^m\to\R$, to be introduced below in the context
of the Lagrange problem),
$\xi^0\in\R^n$ and $U\subseteq\R^m$ are prescribed,
and $\mathcal U$ is the set of piecewise-continuous functions
$u:[0,1]\to U$, 
that have at most a finite number of points
of discontinuity in $[0,1]$ and are right-continuous 
on [0,1).\footnote{
While stricter than necessary, this 
assumption on $\mathcal U$ allows for simpler analysis and is adequate 
for the large majority of applications.}
Of course, given an {\em optimal} (for~(\ref{eq:MayerInt})) control 
$u^*\in\mathcal U$,
the associated optimal state trajectory $x^*(\cdot)$ is 
well defined (and continuous) on [0,1].

Problem~(\ref{eq:MayerInt}) is 
closely related to the {\em finite}-dimensional optimization problem
\begin{equation}
\label{eq:min psi}
\min \psi(\xi) {\rm ~s.t.~} \xi\in K(\xi^0),
\end{equation}
where $K(\xi^0)$ is the set of states reachable at time $t=1$
when starting from state $\xi^0$ at $t=0$ and applying
a control $u\in{\mathcal U}$.  
In particular, if $\xi^*$ solves~(\ref{eq:min psi}),
every control $u^*\in\mathcal U$ that steers $x$
from $x(0)=\xi^0$ to $x(1)=\xi^*$, 
solves~(\ref{eq:MayerInt}).
Now, it is well known (e.g.,~\cite[Chapter V, Lemma 1]{Varaiya72:0})
that if $\xi^*$ solves~(\ref{eq:min psi}) then
\begin{equation}
\label{eq:psi'>=0}
\begin{aligned}
{\rm J}\psi(\xi^*)d&\geq0\\
&~\forall d \in {\rm cl~coCC}(\xi^*,K(\xi^0)),
\end{aligned}
\end{equation}
where ${\rm CC}(\xi,S)$ 
is G.L.~Bouligand's {\em contingent cone}\footnote{\label{foot2}For every 
$S\subset\R^n$ and $\xi\in{\rm cl}(S)$,
${\rm CC}(\xi,S)$ is closed; 
e.g.~\cite[Proposition 6.2]{RW98}.}~\cite[section 3]{Bouligand1928}
to set $S$ at $\xi\in S$.
For ease of reference,
we recall here the definition of the contingent cone, as formulated, e.g.,
in~\cite[Theorem 6.3.6(a)]{BorweinLewis2000}.
\begin{definition}(e.g., \cite[Theorem 6.3.6(a)]{BorweinLewis2000})
\label{def:CC}
Given $S\subset\R^n$ and $\xi\in{\rm cl}(S)$,
the {\em contingent cone}
${\rm CC}(\xi,S)$ consists of those vectors
$d\in\R^n$ for which there are sequences $\{\epsilon_k\}\downarrow0$ in $\R$
and $\{d_k\}\to d$ in $\R^n$ such that $\xi+\epsilon_k d_k\in S$ for all $k$.
\end{definition}

\noindent
Based on a ``needle'' perturbation, sensitivity analysis shows
that, given any $v\in U$, and 
any $\tau\in[0,1)$ at which $u^*$ is continuous, 
\begin{equation}
\label{eq:htauv}
\begin{aligned}
d_{\tau,v}&:=\Phi (1,\tau)\bigl[f(\tau,x^*(\tau),v)-f(\tau,x^*(\tau),u^*(\tau))\bigr]\\
&\in {\rm CC}(x^*(1),K(\xi^0)),
\end{aligned}
\end{equation}
where $x^*$ is the state trajectory associated to $u^*$ 
and $\Phi:[0,1]\times[0,1]\to\R^{n\times n}$ the state-transition 
operator for the linearization of $f$ along $(x^*(t),u^*(t))$;
see, e.g., \cite[Exercise~1 in section VIII.1]{Varaiya72:0}, 
\cite[equation~(4.23)]{Liberzon2012}.   

The following instance of PMP ensues
(see, e.g., \cite[section VIII.1, Theorem 2]{Varaiya72:0}),
where $H:\R\times\R^n\times\R^n\times\R^m\to\R$ 
is the {\em pre-Hamiltonian} 
associated to Problem~(\ref{eq:MayerInt}), given by
\[
H(\tau,\xi,\pi,\upsilon):=f(\tau,\xi,\upsilon)^{\rm T}\pi
\]
and
${\rm J}_\xi f:\R\times\R^n\times\R^m\to\R^n\times\R^n$ denotes
the Jacobian of $f$ with respect to its second argument.
\begin{fact}                                                
\label{thm:FreeTerminalPoint}
Let $u^*\in\mathcal U$ be optimal for problem~(\ref{eq:MayerInt}), 
with associated state trajectory $x^*$.  
Then there exists an
absolutely continuous $p^*:\R\to\R^n$ such that
\[
\begin{aligned}
H(t,x^*(t),p^*(t),u^*(t)) &= 
\max_{\upsilon\in U} H(t,x^*(t),p^*(t),\upsilon)\\
&\forall t\in[0,1),
\end{aligned}
\]
and $p^*$ solves the linear, terminal-value problem
\[
\begin{aligned}
\dot p^*(t) &= -{\rm J}_{\xi}f(t,x^*(t),u^*(t))^{\rm T}p^*(t)\\
&\qquad {\rm a.e.~}t\in[0,1],\quad p^*(1)=-\nabla\psi(x^*(1)).
\end{aligned}
\]

\end{fact}

\noindent
This result is easily extended to problems where $x^*(0)$, rather
than being prescribed, is merely restricted to lie on a specified manifold,
namely, $g^0(x^*(0))=0_{m_0}$, where 
$g^0:\R^n\to\R^{m_0}$, for some $m_0$, is continuously differentiable,
with Jacobian J$g^0:\R^n\to\R^{m_0\times n}$.
For such problem, it is readily
determined 
that, under the condition 
that $x^*(0)$ is regular for $g^0$ (i.e., that
J$g^0(x^*(0))$ is surjective),
Fact~\ref{thm:FreeTerminalPoint} still holds as stated with, 
in addition, $p^*(0)$ being now constrained to satisfy the 
``transversality'' condition\footnote{The $m_0$ degrees 
of freedom ``gained'' by $x^*(0)$ are thus ``lost'' by $p^*(0)$.}
\[
p^*(0)\perp{\mathcal N}\left({\rm J}g^0(x^*(0))\right).
\]

The case of a constrained {\em terminal} state is more challenging.  
Indeed, consider
\begin{equation}
\begin{aligned}
\text{minimize}_{u\in\mathcal U}\quad &\psi(x(1)) \\
\text{s.t.}\quad &\dot x(t)=f(t,x(t),u(t)),\quad {\rm a.e.~}t\in[0,1],\\
&x~{\rm abs.~cont.},\quad g^0(x(0))=0_{m_0},\\
&g^{\rm f}(x(1))=0_{m_{\rm f}}.
\end{aligned}
\label{eq:MayerConstrained}
\end{equation}
where $g^{\rm f}:\R^n\to\R^{m_{\rm f}}$, for some $m_{\rm f}$, 
is continuously differentiable.
In this case, (\ref{eq:psi'>=0}) is not implied anymore,
since not all of the reachable set $K(T_0)$ is now allowable; here 
$K(T_0)$ denotes the set of states reachable at $t=1$ from an 
initial point in 
\[
T_0:=\{\xi\in\R^n:g^0(\xi)=0_{m_0}\}
\]
at $t=0$.
Accordingly,
Fact~\ref{thm:FreeTerminalPoint} is not guaranteed to hold when
the terminal point is constrained.  
It is a well-established fact that the following augmented version
of that result does hold;
see, e.g.,~\cite[section VIII.1, p. 156, Theorem 3]{Varaiya72:0}.


\begin{fact}                                                   
\label{Pontryagin}
Let $u^*\in\mathcal U$ be optimal for problem~(\ref{eq:MayerConstrained}),
with associated state trajectory $x^*$.  
Assume that the initial state $\xi^0~(=x^*(0))$ is a regular point for 
$g^0$ and that ${\rm J}g^{\rm f}(x^*(1))$ is surjective.
Then there exist a constant $p_0^*\geq0$ 
and an absolutely continuous function $p^*:[0,1]\to\R^n$,
with $[p_0^*;p^*(\cdot)]$ not identically zero, such that
\[
\begin{aligned}
H(t,x^*(t),p^*(t),u^*(t))
&=\max_{\upsilon\in U}H(t,x^*(t),p^*(t),\upsilon)\\
&\qquad \forall t\in[0,1),
\end{aligned}
\]
\begin{equation}
\label{eq:costate}
\dot p^*(t) = -{\rm J}_\xi f(t,x^*(t),u^*(t))^{\rm T} p^*(t) 
~{\rm a.e.~} t\in[0,1]
\end{equation}
and the transversality conditions 
\[
\begin{aligned}
p^*(0)&\perp{\mathcal N}\left({\rm J}g^0(x^*(0))\right),\\
p^*(1)+p_0^*\nabla\psi(x^*(1))
&\perp{\mathcal N}\left({\rm J}g^{\rm f}(x^*(1))\right).
\end{aligned}
\]
hold.
\end{fact}
\begin{remark}
Clearly, in the statement of Fact~\ref{Pontryagin}, $p_0^*\in\{0,1\}$ can
be substituted for $p_0^*\geq0$ via scaling of $[p_0^*,p^*]$ by $p_0^*$
when $p_0^*>0$.
\end{remark}

\noindent
These necessary conditions of optimality are ``weak'' 
({\sl a la} FJ Conditions~\cite{john1948extremum}; also see, 
e.g.,~\cite[section 3.3.5]{Bertsekas1999}) in the sense that, 
by not guaranteeing that some {\em nonzero} $p_0^*$ satisfies the conditions,
they leave open the possibility that the conditions may hold regardless of 
the cost function~$\psi$. 

\section{Sufficient Conditions for Normality} 
\label{sec:p0NOT0}

As noted above, 
in the context of Mayer's problem~(\ref{eq:MayerInt}),
in the free-terminal-state case 
(i.e., $g^{\rm f}\equiv0_{m_{\rm f}}$ or, equivalently, $m_{\rm f}=0$),
PMP follows from noting that $u^*\in{\mathcal U}$
is optimal if and only if the associated state trajectory $x^*$ is
such that 
\begin{equation}
\label{eqn:argmin T0}
x^*(1)\in{\rm argmin}_\xi\{\psi(\xi) {\rm ~s.t.~} \xi\in K(T_0)\},
\end{equation}
and a standard first order necessary condition for~(\ref{eqn:argmin T0})
to hold is that
\begin{equation}
\label{eq:CC(K)}
\begin{aligned}
{\rm J}\psi(x^*(1))d&\geq0\\
&\quad\forall d\in{\rm cl~coCC}(x^*(1),K(T_0)).
\end{aligned}
\end{equation}
Clearly, when the terminal state is required to lie within 
\[
T_{\rm f}:=\{\xi\in\R^n:~g^{\rm f}(\xi)=0_{m_{\rm f}}\},
\]
(\ref{eqn:argmin T0}) becomes
\[
x^*(1)\in{\rm argmin}_\xi\{\psi(\xi) {\rm ~s.t.~} \xi\in
T_{\rm f}\cap K(T_0)\},
\]
yielding the necessary condition
\begin{equation}
\label{eq:CC(KcapTf)}
{\rm J}\psi(x^*(1))d \geq 0
\quad\forall d \in {\rm D}:= {\rm cl~coCC}(x^*(1),T_{\rm f}\cap K(T_0)),
\end{equation}
which is {\em weaker} ({\em less} constraining on $x^*$) than~(\ref{eq:CC(K)}), 
since it asserts nonnegativity of the directional derivative over
a {\em subset} of the directions $d$ that satisfy~(\ref{eq:CC(K)}).

While~(\ref{eq:CC(KcapTf)}) allows for $p_0^*=0$, two conditions are now 
introduced which, as shown below, together guarantee the validity of
$p_0^*=1$.


\smallskip\noindent
{\bf Conditions Terminal-State (TS):} 
\begin{enumerate}
\item[TS1.]
${\rm J}\psi(x^*(1))d \geq 0
\quad\forall d \in ~
{\mathcal N}\left({\rm J}g^{\rm f}(x^*(1))\right)
\cap {\rm D}$  
\item[TS2.]
The cone
\[
C := \left\{ \left[
      \begin{array}{c} {\rm J}\psi(x^*(1))  \\
                  {\rm J}g^{\rm f}(x^*(1)) 
    \end{array} \right]{\rm D}
    \right\} \subseteq\R^{1+m_{\rm f}}
\]
is closed.
\end{enumerate}

\noindent
\begin{theorem}
\label{thm:abnormal}
Let $u^*\in\mathcal U$ be optimal for problem~(\ref{eq:MayerConstrained}), 
with corresponding state trajectory $x^*$.   
Suppose $x^*(0)$ is a regular point of $g^0$, ${\rm J}g^{\rm f}(x^*(1))$
is surjective, and Conditions TS hold.  
Then the conditions in Fact~\ref{Pontryagin} hold with $p_0^*=1$.
\end{theorem}

\begin{proof}
Write 
\[
{\rm M}:= \left[
      \begin{array}{c} {\rm J}\psi(x^*(1))  \\
                  {\rm J}g^{\rm f}(x^*(1))
      \end{array} \right],
\quad
{\rm R}:= \left[ \begin{array}{c} -1 \\  0_{m_{\rm f}} \end{array} \right].
\]
Since $x^*(1)$ minimizes $\psi$ over $T_{\rm f}\cap K(T_0)$,
(\ref{eq:CC(KcapTf)}) holds.
Then, from Condition TS1, ${\rm R}\not\in C$;
indeed, if ${\rm R}\in C$, then from the definition of $C$, 
there would exist $d\in{\rm D}$ 
such that J$\psi(x^*(1))d=-1$ and 
J$g^{\rm f}(x^*(1))d=0_{m_{\rm f}}$,
which would contradict TS1; hence, ${\rm }\not\in C$.
Since $C$ is a closed convex cone (TS2), $C$ and $\rm R$ are separated
by a hyperplane through the origin that does not contain $\rm R$, 
i.e., there exist 
$p_0^*\in\{0,1\}$ and $\pi\in{\rm R}^{m_{\rm f}}$ such that,
with $\mu:=[p_0^*;\pi]$, we have
$\mu^{\rm T}{\rm }R<0$ (i.e., $p_0^*=1$) 
and $\mu^{\rm T} w\geq0$ for all
$w\in C$, i.e.,
\[
\left({\rm J}\psi(x^*(1))+\pi^{\rm T}{\rm J}g^{\rm f}(x^*(1))\right)
d \geq 0 \quad\forall d\in {\rm D}.
\]
In particular,
from~(\ref{eq:htauv}), 
\[
\begin{aligned}
&\left({\rm J}\psi(x^*(1))+\pi^{\rm T}{\rm J}g^{\rm f}(x^*(1))\right)\\
&\quad\cdot\Phi(1,\tau)\bigl[f(\tau,x^*(\tau),v)-f(\tau,x^*(\tau),u^*(\tau))\bigr]\\
&\geq0
\end{aligned}
\]
for every $v\in U$, and every $\tau\in[0,1)$ at which $u^*$ is continuous.
Accordingly, 
invoking~\eqref{eq:costate} and
\[
p^*(1) =-\nabla\psi(x^*(1))
-{\rm J}g^{\rm f}(x^*(1))^{\rm T}\pi \quad{\rm for~some}~\pi\in\R^{m_f}
\]
or, equivalently, the {\em transversality} condition
\[
p^*(1)+\nabla\psi(x^*(1))
\perp {\mathcal N}\left({\rm J}g^{\rm f}(x^*(1))\right),
\]
completes the proof.
\end{proof}

\begin{remark}
\label{rem:TS2}
A sufficient condition for TS2 to hold regardless of the
state dynamics is that
{\rm M} be square and invertible, i.e., that $m_{\rm f}=n-1$ and
$\nabla\psi(x^*(1))\not\in{\mathcal R}({\rm J}g^{\rm f}(x^*(1))^{\rm T})$.
Hence if Condition TS1 holds (and ${\rm J}g^{\rm f}(x^*(1))$ is surjective), 
Fact~\ref{Pontryagin} will hold with $p_0^*=1$.
(This is the case for every Lagrange {\em steering} problem, discussed in
sections~\ref{sec:LTIstate} and~\ref{sec:LTIdynANDcost} below.)
\end{remark}

\section{The Lagrange Problem}
\label{sec:Lagrange Problem}

Consider now the {\em Lagrange} problem
\begin{equation}
\begin{aligned}
\text{minimize}_{u\in\mathcal U}\quad
&\int_0^1 f_0(\tau,x(\tau),u(\tau))\,d\tau \\
\text{s.t.}\quad
&\dot x(t)=f(t,x(t),u(t)),\quad {\rm a.e.~}t\in[0,1],\\
&x~{\rm abs.~cont.},\quad g^0(x(0))=0_{m_0},\\
&g^{\rm f}(x(1))=0_{m_{\rm f}}.
\end{aligned}
\label{eq:La1}
\end{equation}
where the closed set $U\subseteq\R^m$ is prescribed, 
$f_0:\R\times\R^n\times\R^m\to\R$ and, as before,
$f:\R\times\R^n\times\R^m\to\R^n$, $g^0:\R^n\to\R^{m_0}$, and
$g^{\rm f}:\R^n\to\R^{m_{\rm f}}$; and
where the same regularity conditions are assumed on $f$, $g^0$, and $g^{\rm f}$
as was done in section~\ref{sec:background}, and 
every condition that is assumed on $f$ is also assumed on~$f_0$.
This problem is readily transposed to the Mayer formalism,\footnote{Similarly,
the Mayer problem is readily transposed to the Lagrange formalism, with cost
integrand given by ${\rm J}\psi(x(\tau))f(\tau,x(\tau),u(\tau))$.}
with augmented state-vector $\tilde x:=[x_0;x_1;\ldots;x_n]$
(with initial state $\tilde x(0)$ set to $[0;x(0)]$, without loss of
generality), yielding, with cost $\tilde\psi(\tilde x(1)):=x_0(1)$,
\begin{equation}
\begin{aligned}
\text{minimize}_{u\in\mathcal U}&\ ~x_0(1) \\
\text{s.t.~}&\ \dot{\tilde x}(t)=\tilde f(t,\tilde x(t),u(t)),{\rm~a.e.~}t\in[0,1],\\
&\ \tilde x~{\rm abs.~cont.},\\
&\ \tilde g^0(\tilde x(0))=0_{1+m_0},\\
&\ \hspace{-1.25em}\tilde g^{\rm f}(\tilde x(1))=0_{1+m_{\rm f}}.
\end{aligned}
\label{eq:La2}
\end{equation}
which involves vector field $\tilde f:\R\times\R^{1+n}\times\R^m\to\R^{1+n}$ defined by 
$\tilde f(\theta,\tilde\xi,\upsilon):=[f_0(\theta,\xi,\upsilon);f(\theta,\xi,\upsilon)]$, 
with associated reachable set 
{\small
\[
\tilde K(\tilde T_0):=\{\tilde\xi'\in\tilde K(\tilde\xi):\tilde\xi\in\tilde T_0\},\,{\rm with}\,
\tilde T_0:=\{\tilde\xi:\tilde g^0(\tilde\xi)=0\}.
\]
}
and mappings $\tilde g^0:\R^{1+n}\to\R^{1+m_0}$ and\linebreak
$\tilde g^{\rm f}:\R^{1+n}\to\R^{m_{\rm f}}$, defined by
\[
\begin{aligned}
\tilde g^0(\tilde\xi)&:=[\xi_0;g^0(\xi)],\\
\tilde g^{\rm f}(\tilde\xi)&:=g^{\rm f}(\xi).
\end{aligned}
\]
Fact~\ref{Pontryagin} applied to this Mayer problem yields, after
simplification, the following result 
(see, e.g., \cite{Varaiya72:0}), which involves the modified pre-Hamiltonian
$\tilde H:\R\times\R^n\times\R^{1+n}\times\R^m\to\R$:
\[
\tilde H(\tau,\xi,\tilde\pi,\upsilon):=
-\pi_0 f_0(\tau,\xi,\upsilon)+\pi^{\rm T}f(\tau,\xi,\upsilon),
\]
with $\tilde\pi:=[\pi_0;\pi]$, $\pi_0\in\{0,1\}$.
\begin{fact}
\label{thm:p0*}
Let $u^*\in\mathcal U$ be optimal for 
problem~(\ref{eq:La1}) with $g^{\rm f}(\xi):=\xi-\xi^{\rm f}$,
with associated state trajectory
$x^*$.  Then, if
${\rm J}g^{\rm f}(x^*(1))$ is surjective,
there exists an absolutely continuous $p^*:\R\to\R^n$
such that
\[
\begin{aligned}
\tilde H(t,x^*(t),\tilde p^*(t),u^*(t))
&=\max_{v\in U} \tilde H(t,x^*(t),\tilde p^*(t),v)\\
&~\forall t\in[0,1)
\end{aligned}
\]
with constant $p_0^*\in\{0,1\}$, $\tilde p^*:=[p_0^*;p^*]$ not identically zero
and satisfying the linear equation
\[
\begin{aligned}
\dot p^*(t)&=-{\rm J}_\xi f(t,x^*(t),u^*(t))^{\rm T}p^*(t)\\
&\quad+p_0^*\nabla_\xi f_0(t,x^*(t),u^*(t))\\
&\qquad {\rm a.e.~}t\in[0,1).
\end{aligned}
\]
as well as the transversality conditions
\[
\begin{aligned}
p^*(0)&\perp{\mathcal N}\left({\rm J}g^0(x^*(0))\right),\\
p^*(1) &\perp{\mathcal N}\left({\rm J}g^{\rm f}(x^*(1))\right).
\end{aligned}
\]
\end{fact}

\begin{remark}
Of course, in the case of a Lagrange problem, 
since the cost function in the associated Mayer problem is $x_0(1)$,
the terminal value $\tilde\xi:=\tilde x^*(1)$ of the optimal augmented-state 
trajectory 
must lie
on the ``lower boundary'' of $\tilde K(\tilde T_0)$, i.e.,
$\tilde\xi^*=[\xi^*_0;\xi^*]$
must satisfy $\xi_0^*=\min\{\xi'_0:[\xi_0';\xi^*]\in\tilde K(\tilde T_0)\}$.
\end{remark}

With $\tilde T_0:=\{[0;\xi^0]\}=\{\tilde \xi^0\}$, we now
simplify notation and define
\[
\begin{aligned}
K(t;\xi^0)&:=K(t;T_0), \quad \tilde K(t;\xi^0):=\tilde K(t;\tilde T_0);\\
K(\xi^0)&:=K(1;\xi^0), \quad \tilde K(\xi^0):=\tilde K(1;\xi^0).
\end{aligned}
\]
We also define $\tilde T_{\rm f}:=\{\tilde \xi\in\R^{1+n}:g^{\rm f}(\xi)=0\}$

A simpler rewrite of Condition TS1
for the case of (the Mayer reformulation of) a Lagrange problem
is as follows.
\begin{enumerate}
\item[LTS1.] cl co${\rm CC}(\tilde x^*(1),\tilde K(T_0))
\cap{\mathcal N}({\rm J}\tilde g^{\rm f}(\tilde x^*(1)))
\cap\{[d_0;d]\in\R^{1+n}:d_0<0\}=\emptyset.$
\end{enumerate}
Geometrically, Condition LTS1 says that the tangent plane to
$\tilde T_{\rm f}$ at $\tilde x^*(1)$ contains no ``downward tangent'' 
(vertical or otherwise) to $\tilde K(T_0)$ at $\tilde x^*(1)$.
\begin{remark}
The above points to a well-known close connection between 
normality of the PMP and finite sensitivity
of the optimal cost to small changes in the problem data (see, e.g.,
\cite{Babad1991,BabadVinter1993}.)  Indeed, it is
clear that, except for very special cases, a slight perturbation
of the data can result in $\tilde T_{\rm f}$ either ``moving out'' of
$\tilde K(T_0)$, thus making the ``target'' $\tilde T_{\rm f}$ locally
unreachable, or ``moving into''
$\tilde K(T_0)$, causing its ``lowest'' intersection with $\tilde K(T_0)$,
whose $0$th coordinate is the optimal cost $x_0^*(1)$, to make a quantum
``jump'', indicating infinite sensitivity.
\end{remark}

Next, when Condition TS2 is applied to the Lagrange problem, 
since ${\rm J}\tilde g^{\rm f}(x^*(1))=[0~~{\rm J}g^{\rm f}(x^*(1))]$,
matrix $\rm M$ defined at the top of the proof of Theorem~\ref{thm:abnormal}
reads
\begin{equation}
\label{eq:M}
{\rm M}
= \begin{bmatrix} 1 & 0_n \\ 0 & {\rm J}g^{\rm f}(x^*(1)) \end{bmatrix},
\end{equation}
so that the condition is satisfied if and only if, beside being surjective (as
assumed in Fact~\ref{thm:p0*}), ${\rm J}g^{\rm f}(x^*(1))$ also is 
injective, i.e., it is invertible.
In particular, for such ``steering'' Lagrange problem, it follows from
Theorem~\ref{thm:abnormal} that
Fact~\ref{thm:p0*} holds with $p_0^*=1$ subject only to Condition LTS1.

The next theorem plays a central role in the remainder of this paper.
It involves the ``lower-boundary'' function $\nu:K(T_0)\to\R$ defined by 
\[
\nu(\xi):=\inf\{\xi_0:[\xi_0;\xi]\in\tilde K(T_0)\}.
\]

\begin{theorem}
\label{thm:LTS1whenri(K)}
Let $\tilde x^*:[0,1]\to\R^{1+n}$ be an optimal augmented-state trajectory.
Suppose that (i) $x^*(1)\in{\rm relint}(K(T_0))$,
(ii) $\nu$ is convex in a neighborhood of $x^*(1)$, and 
(iii) ${\rm J}g^{\rm f}(x^*(1))$
and ${\rm J}g^{\rm f}|_{{\rm aff}(K(T_0))}$ are surjective.
Then LTS1 holds at $\tilde x^*(1)$.  
\end{theorem}

\begin{proof}
Without loss of generality, reduce $\R^n$ to aff$(K(T_0))$, so that
${\rm relint}(K(T_0))$ becomes the interior of $K(T_0)$.
Since $x^*(1)$ is an interior point of the reduced reachable set,
${\rm CC}(x^*(1),K(T_0))=\R^n$, so every $d\in{\rm D}$ belongs to the
contingent cone.
Since every convex function on an open set is locally 
Lipschitz~\cite[Corollary 2.27]{MR4633252}, it follows from (i) and (ii) that
there exists $\delta>0$ such that, in the reduced
space,  $\nu$ has a Lipschitz constant $L\in\R$ on $B_\delta(x^*(1))$.
Now, proceeding by contradiction, suppose that LTS1 does not hold, i.e.,
there exists $\tilde d=[d_0;d]\in\R^{1+n}$ such that
\begin{itemize}
\item[(a)] $d_0<0$;
\item[(b)] $\tilde d\in{\mathcal N}({\rm J}\tilde g^{\rm f}(\tilde x^*(1)))$,
i.e., $d\in{\mathcal N}({\rm J}g^{\rm f}(x^*(1)))$; and
\item[(c)] $\tilde d\in{\rm CC}(\tilde x^*(1),\tilde K(T_0))$.
\end{itemize}
By (b) and (iii), there exists $\Gamma:\R\to T_{\rm f}\cap B_\delta(x^*(1))$,
continuously differentiable,
such that $\Gamma(0)=x^*(1)$ and $\Gamma'(0)=d$.  By (c),
there exists $\{\epsilon^k\}\subset\R$ and $\{\tilde o^k\}\subset\R^{1+n}$
such that, as $k\to\infty$,
$\epsilon^k\downarrow0$ and $\frac{\|\tilde o^k\|}{\epsilon^k}\to0$
and, for all $k$, 
$\tilde\xi^k := \tilde x^*(1)+\epsilon^k\tilde d +\tilde o^k\in\tilde K(T_0)$
and $\xi^k\in B_\delta(x^*(1))$.
Since, for all $k$, $\nu(\Gamma(\epsilon^k))\geq\nu(x^*(1))$ 
(by optimality of $\tilde x^*(\cdot)$) and $\nu(\xi^k)\leq\xi_0^k$ 
(by definition of $\nu$), the Lipschitz constant $L$ satisfies
\[
L\geq\frac{\nu(\Gamma(\epsilon^k))-\nu(\xi^k)}{\|\Gamma(\epsilon^k)-\xi^k\|} 
\geq\frac{\nu(x^*(1))-\xi_0^k}{\|\Gamma(\epsilon^k)-x^*(1)-\epsilon^k d 
- \tilde o^k\|}
\]
implying, from the definition of $\tilde\xi^k$,
\[
L\geq\frac{\nu(x^*(1)) - x_0^*(1) - \epsilon^k d_0 - o_0^k}
{\|\Gamma(\epsilon^k)-x^*(1)-\epsilon^k d - o^k\|}
\]
and since $\nu(x^*(1))=x_0^*(1)$, dividing up and down by $\epsilon^k$,
\[
L\geq\frac{-d_0-\frac{o_0^k}{\epsilon^k}}
{\left\|
\frac{\Gamma(\epsilon^k)-x^*(1)}{\epsilon^k}- d - \frac{o^k}{\epsilon^k}
\right\|} =:E_k ~~\forall k.
\]
Since (by (a)) $-d_0>0$ and since the denominator of $E_k$ tends to $0+$
(since $d=\Gamma'(0)$),
$E_k\to+\infty$ as $k\to\infty$, in contradiction to $L\in\R$ being a
(finite) Lipschitz constant.
\end{proof}

In the next three subsections, still in the context of Lagrange problems,
we successively narrow down the scope of the analysis to 
(i) problems with {\em linear state dynamics}, 
(ii) {\em steering} (to a single point) problems with {\em LTI state dynamics}, 
and finally 
(iii) steering problems
with LTI state dynamics {\em and LTI cost integrand} (in the sense that
$f(\tau,\xi,\upsilon)$ does not vary with $\tau$ and is linear in 
$(\xi,\upsilon)$.  Further, we focus on the case of scalar $u$, and
make the following
\begin{assumption}
From this point onward, the proof is specialized to the normalized symmetric
control interval $U=[-1,1]$.  This specialization is essential for the
explicit formulas $u={\rm sgn}(\cdot)$ and for the perturbation identity
$u^k-u^*\in\{0.\pm2\}$ used below to hold.
\end{assumption}


\subsection{Lagrange Problem with Linear State Dynamics}
\label{sec:LTI dynamics}

Here, we focus on vector fields $f$ of the form
\begin{equation}
\label{eq:LTI dynamics}
f(\tau,\xi,\upsilon):=A(\tau)\xi+b(\tau)\upsilon
\end{equation}
for some continuous mappings $A:\R\to\R^{n\times n}$ and $b:\R\to\R^n$, and
we let $\Phi_A$ be the state-transition matrix associated with $A$.
We assume, without loss of generality, that $x^*(0)=0_n$.  
There indeed is no loss of generality 
in doing so since, 
with the change of variable $x^*\mapsto y^*$, with
$y^*:=x^*-x^*(0)$ (so that $y^*(0)=0_n$), $y^*-x^*$ is independent of $u$
(hence the change of variable amounts to a mere translation of 
$K(x^*(0))$): we have
\[
\begin{aligned}
\dot y^*(t)=\dot x^*(t) &=  A(t)x^*(t) + b(t)u(t) \\
&= A(t)y^*(t) + A(t)x^*(0)+ b(t)u(t),
\end{aligned}
\]
yielding, since $y^*(0)=0_n$,
\[
\begin{aligned}
y^*(t)&=\int_0^t\Phi_A(t,\tau)
\bigl(A(\tau)x^*(0)+b(\tau)u(\tau)\bigr)d\tau\\
&=\int_0^t\Phi_A(t,\tau)b(\tau)u(\tau)d\tau\\
&\quad+\int_0^t\Phi_A(t,\tau)A(\tau)x^*(0)d\tau,
\end{aligned}
\]
and, since $x^*(t)=\Phi_A(t,0)x^*(0)+\int_0^t\Phi_A(t,\tau)b(\tau)u(\tau)d\tau$,
we get
\[
\begin{aligned}
y^*(t) - x^*(t) &= \Phi_A(t,0)x^*(0)\\
&\quad+\int_0^t\Phi_A(t,\tau)A(\tau)x^*(0) d\tau,
\end{aligned}
\]
which is independent of $u$, as claimed.

Accordingly, for the sake of notational simplicity, we now define
\[
K:=K(0_n),
\quad\tilde K:=\tilde K(0_n).\]
Finally, given $u\in\mathcal U$, we denote by $\tilde x(t;u)$
the augmented state at time $t$, starting from $\tilde x(0)=0_{1+n}$ at
time $t=0$, and applying control $u$.  In particular, for $u\in\mathcal U$,
\begin{equation}
\label{eq:x(1;u)}
x(1;u) = \int_0^1 \varphi(t) u(t) dt,
\end{equation}
where $\varphi:\R\to\R^n$ is defined by
\[
\varphi(t):=\Phi_A(1,t)b(t).
\]

Now, as an application of Theorem~\ref{thm:LTS1whenri(K)}
consider the case with $f$ given by~(\ref{eq:LTI dynamics}) and
\begin{equation}
\label{eq:f0=ax+h}
f_0(\tau,\xi,\upsilon):=a_0(\tau)^{\rm T}\xi+h_0(\tau,\upsilon),
\end{equation}
where $h_0:\R\times\R\to\R$ and $a_0:\R\to\R^n$ are continuous
(as required by our assumption on $f_0$).

\noindent
The following 1963 result is due to
L.W.~Neustadt~\cite[Theorem 1]{Neustadt1963},\footnote{Neustadt's
result also applies to the case of state dynamics more general 
than~(\ref{eq:LTI dynamics}):
$f(\tau,\xi,\upsilon)=A(\tau)\xi + h(\tau,\upsilon)$,
with $h:\R\times\R\to\R^n$ continuous.}
where a 1951 result by D.~Blackwell~\cite{Blackwell1951} plays a major role;
also see~\cite[Theorem 1]{Pecsvaradi1971}.
\begin{lemma}
\label{lem:Neustadt}
Given dynamics~(\ref{eq:LTI dynamics}) and cost integrand~(\ref{eq:f0=ax+h})
with $h:\R\times \R\to\R$ continuous,
if $U$ is compact, then $K$ and $\tilde K$ are convex and compact.
\end{lemma}

\noindent
It follows that, under the same assumptions, $\nu$ is convex
and hence Theorem~\ref{thm:LTS1whenri(K)} applies, i.e., LTS1 holds
whenever $x^*(1)\in{\rm relint} K(T_0)$ and ${\rm J}g^{\rm f}(x^*(1))$
is surjective.  

\subsection{Lagrange Steering Problem with LTI state dynamics}
\label{sec:LTIstate}

We now focus on the class of Lagrange problems 
with $g^{\rm f}(\xi):=\xi-\xi^{\rm f}$ 
(and w.l.o.g. $x^*(0)=0$: see section~\ref{sec:LTI dynamics}), 
where $\xi^{\rm f}\in K$ is prescribed, 
which we refer to as the Lagrange {\em steering} problem.
In such case, ${\rm J}g^{\rm f}(\xi)=I_n$, ${\rm M}=I_{1+n}$
(see~(\ref{eq:M})), and TS2 is always satisfied (see Remark~\ref{rem:TS2}).
From here on, time-invariance is also assumed, i.e., the
state dynamics are now
\begin{equation}
\label{Ax+Bu}
f(\xi,\upsilon)=A\xi + b\upsilon.
\end{equation}
In such a case, the overall picture simplifies, as noted next.

\begin{lemma}
\label{lem:reduced(A,b)}
A (necessary and) sufficient condition for {\rm LTS1} to hold for a
given pair $(A,b)$ is that these conditions hold for the associated 
controllable subsystem $(A',b')$, 
with $A'$ in real Jordan canonical form.
\end{lemma}
\begin{proof}
First note that, if $(A,b)$ is {\em not} controllable, then
$\tilde K$ entirely lies in a proper subspace of $\R^{1+n}$ whose ``vertical''
projection is the controllable subspace of $\dot x=Ax+bu$.  An associated
controllable reduced version of $(A,b)$ is then readily obtained as
the controllable portion $(A_{\rm c},b_{\rm c})$ of the Kalman decomposition 
of $(A,b)$ (obtained via a change of basis).  
Next, via an additional change of basis (which of course also
transforms $b$), $A$ can be brought to real Jordan canonical form.  
(In such form, diagonal blocks 
associated with real eigenvalues are real scalars (the real eigenvalues), 
while for non-real eigenvalues 
$\alpha+{\rm j}\omega$ (with $\alpha$ and $\omega$ real, $\omega\not=0$) 
they take the form $[\alpha~~\omega;-\omega~~\alpha]$;
off-diagonal blocks are all zero except for occasional super-diagonal
1s squeezed between the scalar blocks associated to 
two equal real eigenvalues, and occasional 
super-diagonal $I_2$s squeezed between two $2\times 2$ diagonal blocks
associated to equal non-real eigenvalues.)  
Since the uncontrollable portion does not affect the reachable set $K$, 
nor indeed does it affect whether or not LTS1 holds, the proof is complete.
\end{proof}

\noindent
Note that, with $A'$ in real Jordan canonical form, controllability
of $(A',b')$ is equivalent to the pair of conditions: (i) $b_i\not=0$ 
when row $i$ of $A$ is the last row of a Jordan block associated with 
a real eigenvalue and (ii) $[b_i;b_{i-1}]\not=0_2$ when row $i$ of $A$ 
is the last row of a Jordan block associated with a non-real eigenvalue.
(Indeed, the pair $([\alpha~~\omega;-\omega~~\alpha],[b;b'])$,
with $\omega\not=0$, is controllable if and only if $[b;b']\not=0_2.$)

Next, recall that the concept of bang--bang control plays a key role in 
optimal control problems with LTI state dynamics.  In the present
context ($u(t)\in\R~\forall t$, $U$ a compact interval of $\R$) the following
definition is appropriate.
\begin{definition}
A control signal $u\in\mathcal U$ is said to be {\em bang--bang} if
$u(t)\in\{U_{\rm m},U_{\rm M}\}\forall t\in(0,1]$.
\end{definition}
Part of the following result is invoked in the proof of 
Theorem~\ref{col:xfin K(T_0)} below, and all of it in that
of Theorem~\ref{prop:vertical tangent} in the next subsection.
\begin{lemma}
\label{LemmaA}
Consider the (controllable, 
as per Lemma~\ref{lem:reduced(A,b)}) 
Lagrange steering problem with LTI state dynamics, with all
eigenvalues of $A$ in~(\ref{Ax+Bu}) being real, 
and let $U:=[u_{\rm m},u_{\rm M}]$.
Then, given any $\xi\in\partial K$, there exists a (a.e.) unique
$u_{\xi}\in\mathcal U$ such that $x(1;u_{\xi})=\xi$.
Further, $u_{\xi}$ is bang--bang with at most $n-1$ sign changes.
\end{lemma}
\begin{proof}
Let $\xi\in\partial K$ and let
$\gamma\in\R^n\setminus\{0\}$ be an outward normal to a 
supporting hyperplane to $K$ at $\xi$ (Lemma~\ref{lem:Neustadt}).
A control $u_{\xi}\in\mathcal U$
is such that $ x(1;u_{\xi})=\xi$ if and only if
\[
\begin{aligned}
\gamma^{\rm T}\xi&=\max_{\eta\in K}\gamma^{\rm T}\eta\\
&=\max_{u\in\mathcal U}\gamma^{\rm T} x(1;u)\\
&= \gamma^{\rm T} x(1;u_{\xi}),
\end{aligned}
\]
i.e., if and only if
\[
\begin{aligned}
u_{\xi}&\in{\rm argmax}_{u\in\mathcal U} \gamma^{\rm T} x(1;u)\\
&= {\rm argmax}_{u\in\mathcal U} \int_0^1\gamma^{\rm T}\varphi(t)u(t)dt
\end{aligned}
\]
i.e., if and only if $u_{\xi}$ satisfies
\[
u_{\xi}(t) = {\rm sgn}(\gamma^{\rm T}\varphi(t))~~{\rm a.e.~}t\in[0,1].
\]
Now,
$\gamma^{\rm T}\varphi~(=\gamma^{\rm T} e^{(1-\cdot) A} b$),
as an exponential polynomial of order $n$ with real exponents, 
vanishes no more than
$n-1$ times in $\R$ (see, e.g., \cite{AbsilTits2024}) and in particular 
no more than $n-1$ times in [0,1].  The claim follows.
\end{proof}

Now, let $u^*\in\mathcal U$ solve~(\ref{eq:La1})
with $g^{\rm f}(\xi):=\xi-x^{\rm f}$ 
and let $\tilde x^*$ be the associated augmented-state trajectory.
We have
\[
{\mathcal N}({\rm J}\tilde g^{\rm f}(\tilde x^*(1)))= \{[d_0;0_n]:d_0\in\R\},
\]
and Condition LTS1 reduces to
\[
{\rm CC}(\tilde x^*(1),\tilde K)
\cap[d_0;0_n] 
\cap\{[d_0;d]\in\R^{1+n}:d_0<0\}
=\emptyset
\]
so that LTS1 holds at $\tilde x^*(1)$ if and only if every {\em vertical} 
``tangent'' $\tilde d\in{\rm CC}(\tilde x^*(1),\tilde K)$ is upward.
In view of Theorem~\ref{thm:abnormal} and of 
Lemma~\ref{LemmaA},
the hypotheses in 
Theorem~\ref{thm:LTS1whenri(K)} reduce to $\xi^{\rm f}\in{\rm int}(K)$ and
convexity of $\nu$.
We obtain the following.
\begin{theorem}
\label{col:xfin K(T_0)}
Let $u^*\in\mathcal U$ be optimal for Lagrange steering
problem~(\ref{eq:La1}) with $g^{\rm f}(\xi):=\xi-\xi^{\rm f}$,
with LTI state dynamics $f(\xi,\upsilon)=A\xi+b\upsilon$, 
and with cost integrand 
$f_0(\tau,\xi,\upsilon)=a_0^{\rm T}\xi + h(\upsilon)$
i.e.,
\[
\min_{u\in\mathcal{U}}\int_0^1\left(a_0^T x(\tau)+h(u(\tau))\right)d\tau
\]
where $h$ satisfies the assumptions
in Lemma~\ref{lem:Neustadt}.
Then, either (i) $u^*$ is bang--bang or (ii)
Fact~\ref{thm:p0*} holds with $p_0^*=1$ (or both). 
\end{theorem}
\begin{proof}
If $\xi^{\rm f}\in\partial K$, it follows from Lemma~\ref{LemmaA} that $u^*$ is
bang--bang, while if $\xi^{\rm f}\in{\rm relint}(K)$, it follows from 
Remark~\ref{rem:TS2} and Theorem~\ref{thm:LTS1whenri(K)} that 
Fact~\ref{thm:p0*} holds with $p_0^*=1$.
\end{proof}

\noindent
The following simple minimum-energy problem of dimension $n=2$ 
shows that there are cases for which option (ii) does not hold (i.e.,
option (i) cannot be taken out):
\[
\begin{aligned}
{\rm minimize}_{u\in\mathcal U} \int_0^1 
\frac{1}{2} u(t)^2  \,dt 
\end{aligned}
\]
\[
\begin{aligned}
{\rm ~subject~to~}\quad
\dot x &= {\rm diag}(\alpha_1,\alpha_2)x + bu,\\
~x(0)&=0_2, ~x(1;u)=x^{\rm f},\\
~u&\in\mathcal U.
\end{aligned}
\]
Let $K$ be the associated reachable set at time $t=1$ and
let $\xi^{\rm f}\in\partial K$.  It follows from Lemma~\ref{LemmaA}
that every $u\in\mathcal U$ such that $x(1;u)=\xi^{\rm f}$ is bang--bang.  
In particular, every optimal control is bang--bang.
We now show that $p_0^*=0$ indeed.
Pre-Hamiltonian $\tilde H$ reads (see Fact~\ref{thm:p0*}):
\[
\tilde H(\xi,\tilde\pi,\upsilon)=
-\frac{1}{2}\pi_0\upsilon^2+\pi_1b_1\upsilon +\pi_2b_2\upsilon
\]
and setting its $\upsilon$ derivative to zero yields
$
\pi_0\upsilon=\pi_1 b_1 + \pi_2 b_2.
$
Suppose now that 
$p_0^*=1$.  Then, if $u^*$ is an optimal control, we must have
\[
u^*(t)={\rm sat}_{[-1,1]} 
(b_1p_1^*(t)+b_2p_2^*(t))
\]
and $p^*$ must satisfy
\[
\dot p_i^*(t) = -\alpha_i p_i^*(t) \quad i=1,2,
\]
yielding (when $\alpha_1\not=0\not=\alpha_2$)
\[
u^*(t)={\rm sat}_{[-1,1]} 
\left(\sum_{i=1}^2 b_ip_i^*(0)e^{-\alpha_i t}\right)
\]
Clearly, for any given unequal nonzero $\alpha_i$s and any $a_{0i}$, $i=1,2$,
and any $p^*(0)\not=0_n$,
there exists $b\in\R^2$ such that 
the (continuous) argument of ``sat'' changes sign
at some $t\in(0,1)$, resulting in $u^*$ not being bang--bang.  (It is 
readily checked that the same conclusion holds when $\alpha_1\alpha_2=0$.)
With such $b$, Fact~\ref{thm:p0*} holds only with $p_0^*=0$, as claimed.

In view of this negative result, in the remainder of this paper
we limit the study of the case $\xi^{\rm f}\in\partial K$ to Lagrange 
steering problems with LTI dynamics {\em and} LTI cost integrand.

\subsection{Lagrange steering problem with LTI state dynamics,
LTI cost integrand}
\label{sec:LTIdynANDcost}

It is known (e.g.,~\cite[p.258, l. --17 to --15]{BettiolFrankowska2007}) that, 
in the absence of pointwise constraints on the state trajectory (no such
constraints are present in the present paper), 
under ``calmness'' conditions (see, e.g.,~\cite[3rd equation in section 4.2; 
section 6.4]{Clarke1983}),
the Pontryagin Principle is normal indeed (i.e., $p^*_0\not=0$).
In this section, after preliminary results are established, 
a direct proof is given that, in the present context,
regardless of calmness considerations, 
when all coefficients are real as are the eigenvalues of the 
state-dynamics matrix $A$, normality holds; further it still holds 
in some cases with non-real eigenvalues (e.g., in the $2\times 2$ case
when eigenvalues have ``small enough'' imaginary parts).  
See sections~\ref{AwithRealEvalues} 
and~\ref{SomeComplexEvalues}.

Thus, consider a Lagrange steering
problem where both state dynamics and cost integrand are LTI,
viz.,
\begin{equation}
\label{eq:f,f0}
f(\xi,\upsilon)=A\xi+b\upsilon,
\qquad f_0(\xi,\upsilon)=a_0^{\rm T} \xi + b_0\upsilon, 
\end{equation}
where $A$ is square with real entries,
$a_0\in\R^n$ and $b_0\in\R$. For such problem, the restriction 
$\xi^{\rm f}\in{\rm int}(K)$ 
in the statement of Theorem~\ref{thm:LTS1whenri(K)} can be relaxed to
$\xi^{\rm f}\in K$) in many important cases.

Specifically, since $\xi^{\rm f}\in K$, we now 
explore the case when $\xi^{\rm f}\in\partial K$.
(In view of Lemma~\ref{lem:Neustadt}, $\partial K\subseteq K$.)
Given that $\tilde x(0)=0_{1+n}$,
\[
\tilde x(1;u) =\int_0^1 \tilde\varphi(\tau)u(\tau)d\tau,
\quad {\rm with~} \tilde\varphi(t):=e^{(1-t)\tilde A}\tilde b,
\]
where $\tilde A$ is the augmented-state-dynamics matrix.
The next lemma is identical to Lemma~\ref{LemmaA}, with state space
$\R^{1+n}$ 
(including cost integrand) instead of $\R^n$ (see inserted 
tilde characters).
\begin{lemma}
\label{LemmaB}
Assuming that, in addition to $(A,b)$ being controllable
(no loss of generality given Lemma~\ref{lem:reduced(A,b)}),
the augmented system~(\ref{eq:f,f0}) is controllable and all
eigenvalues of $A$ are real, consider the Lagrange steering 
problem with LTI state dynamics
and LTI cost integrand.  
Then, given any $\tilde\xi\in\partial\tilde K$, there exists a unique
(up to equality a.e.)
$u_{\tilde\xi}\in\mathcal U$ such that 
$\tilde x(1;u_{\tilde\xi})=\tilde\xi$.
Further, $u_{\tilde\xi}$ is bang--bang with at most $n$ sign changes.
\end{lemma}

\begin{theorem}
\label{prop:vertical tangent}
With the same assumptions as in Lemma~\ref{LemmaB},
consider the (controllable, as per Lemma~\ref{lem:reduced(A,b)}) 
Lagrange steering problem with LTI state
dynamics, LTI cost integrand, and $\xi^{\rm f}\in\partial K$. 
Let $u^*\in\mathcal U$ be optimal, hence (see Lemma~\ref{LemmaB})
bang--bang with $\bar\imath^*\leq n$ sign changes, and let 
$\{t^*_1,\cdots t^*_{\bar\imath^*}\}\subset[0,1]$
be the times at which $u^*$ changes sign.
Suppose that, given (i) any
$I, I^c\subset\N$, finite index sets, with $I\cap I^c=\emptyset$,
$I\cup I^c\not=\emptyset$, and $\#I^c + 3 \#I \leq n$, and (ii) any
$\{t_i\}_{i\in I}\subseteq\{t_1^*,\ldots,t^*_{\bar\imath^*}\}$ 
and
$\{t_i\}_{i\in I^c}\subseteq[0,1]\setminus\{t_i\}_{i\in I}$, with
$\{t_1^*,\ldots,t^*_{\bar\imath^*}\}\subseteq\{t_i\}_{i\in I\cup I^c}$,
it holds that the matrix
\[
\Phi:= \left[
(\varphi(t_i))_{i\in I^c} ~ ({\rm J}\varphi(t_i))_{i\in I}  
~ (\varphi(t_i))_{i\in I}
\right]
\]
has full column rank.  Then Condition LTS1 holds.
\end{theorem}

\begin{proof}
Let $\tilde x^*$ be the augmented-state trajectory associated to $u^*$.
Proceeding by contraposition, suppose that LTS1 does not hold. 
We proceed to construct an instance where conditions (i) and (ii) are
satisfied but $\Phi$ does not have full column rank, completing the proof.

Given that LTS1 does not hold, since 
${\mathcal N}$(J$\tilde g^{\rm f}(\tilde x^*(1))
=\{\tilde\xi\in\R^{1+n}:\xi=0_n\}$, 
it must hold that $[-1;0_n]\in{\rm CC}(\tilde x^*(1),\tilde K)$, i.e.,
there must exist $\{\epsilon^k\}_{k=1}^\infty\subset\R$ and 
$\{\tilde o^k\}_{k=1}^\infty\subset\R^{1+n}$, $k\in \N$ such that, 
as $k\to\infty$,
$\epsilon^k\downarrow0$, 
$\frac{\|\tilde o^k\|}{\epsilon^k}\to0$ and, for all~$k$,
\[
\tilde\xi^k:=\tilde x^*(1) + \epsilon^k\tilde d + \tilde o^k\in\tilde K
\]
where 
$\tilde d:=[-1;0_n]$.

We now ``groom'' sequence $\{\tilde\xi^k\}$ by the following eight steps.

\smallskip\noindent
{\em Step 1.}
Remove all $\tilde\xi^k$s such that $\xi_0^k\geq x_0^*(1)$.  Since $d_0=-1$,
there are at most finitely many such points.

\smallskip\noindent
{\em Step 2.}
For all $k$ such that $\tilde\xi^k\not\in\partial\tilde K$, move $\tilde\xi^k$
along the direction toward $\tilde x^*(1)+\epsilon^k\tilde d$ until reaching 
$\partial K$.  Since $\tilde x^*(1)+\epsilon^k\tilde d\not\in\tilde K$ (by 
optimality of $\tilde x^*(1)$ and since $d_0=-1$ and $\epsilon^k>0$), 
$\partial K$ is reached before 
$\tilde x^*(1)+\epsilon^k\tilde d$.  Such move does not increase 
$\|\tilde o^k\|$ (since, after the move, the new $\tilde\xi^k$
lies on the line segment from the old $\tilde\xi^k$ ($\in K$)
to $\tilde x^*(1)+\epsilon^k\tilde d$ ($\not\in\tilde K$ since $\tilde x^*(1)$ 
is optimal and $d_0<0$), specifically, at the intersection
of that former line segment with $\partial\tilde K$), hence 
$\frac{\|\tilde o^k\|}{\epsilon^k}\to0$ is preserved.  
Since $K$ is closed (Lemma~\ref{lem:Neustadt}) 
$\tilde\xi^k\in\partial\tilde K$ still holds.  
In view of Step~1 and since $d_0<0$, $\xi_0^k<x_0^*(1)$ still holds for all $k$.

\smallskip\noindent
{\em Step 3.}
By Lemma~\ref{LemmaB}, there is a (a.e.) unique 
control $u^k$ that such that the associated
trajectory satisfies $\tilde x(1;u^k)=\tilde\xi^k$, and $u^k$ is bang--bang with
at most $n$ sign changes.  Let 
$\{\hat t_j^k\}_{j=1}^{\bar\jmath_k}$ be the successive
times of sign changes.  Reduce $\{\tilde\xi^k\}$ to an infinite subsequence
on which $\bar\jmath_k$ is a constant $\bar\jmath$ (such subsequence exists
since $\bar\jmath_k\leq n$ for all $k$) and ${\rm sgn}(u^k(0))$ is a constant
$\hat\sigma$.  Observe that there must exist a nonzero such $\bar\jmath$
since otherwise $\bar\jmath_k=0$ for all $k$ large enough, i.e., for every $k$
large enough, $u^k\in\{\pm1\}$ does not vary with $k$ and hence
$\tilde\xi^k=\tilde x^*(1)$ for all $k$,
in contradiction with 
$\xi_0^k<x_0^*(1)$ for all $k$ (Step 2).  

\smallskip\noindent
{\em Step 4.}
Further reduce $\{\tilde\xi^k\}$ to an infinite subsequence such that 
the associated subsequence $(\hat t_j^k)_{k\in\N}$ converges for
all $j\in\{1,\ldots,\bar\jmath\}$.  
Let 
$\hat t_j^\infty:=\lim_{k\to\infty}\hat t_j^k\in[0,1],~j=1,\ldots,\bar\jmath$
and let $\{t_1,\ldots,t_{\bar\imath}\}:=\{\hat t_j^\infty\}_{j=1}^{\bar\jmath}$,
with $t_1<\cdots <t_{\bar\imath}$.  Further let $i_j$ be defined by
$t_{i_j}=\hat t^\infty_j$ and, for all $t\in[0,1]$, let 
\[
m_t:=\#\{j:\hat t^\infty_j=t\}.
\]
  Also, given any nondecreasing finite-length sequence
$\{\breve t_1,\ldots,\breve t_{\bar m}\}\subset[0,1]$ 
and $\breve\sigma\in\{-1,1\}$, let
\[
u_{\rm BB}(t;\breve\sigma,\{\breve t_n\}_{n=1}^{\bar m})
:= \breve\sigma (-1)^{\#\{m:\breve t_m\leq t\}}, ~\forall t\in[0,1],
\]
a bang--bang signal.
 that, for all $k$,
$u^k\equiv u_{\rm BB}(t;\hat\sigma,\hat t_1^k,\ldots,\hat t_{\jmath}^k)$ and
observe that $u_{\rm BB}(t;\breve\sigma,\{\breve t_n\}_{n=1}^{\bar m})$ is
invariant under the following operations:
\begin{itemize}
\item (i) removal/addition of an even number of repeated $\breve t_m$s;
\item (ii) removal/addition of any number of $\breve t_m$s at $t=1$;
\item (iii) removal of an odd number of repeated $\breve t_m$s at
$t=0$ with concurrent sign change (multiplication by $-1$) of $\breve\sigma$;
\end{itemize}

We now have, by continuity,
\[
\begin{aligned}
&\tilde x\bigl(1;u_{\rm BB}
(\cdot;\hat\sigma,\hat t_1^\infty,\ldots,\hat t_{\bar\jmath}^\infty)\bigr)\\
&\quad=\lim_{k\to\infty}\tilde x\bigl(1;u_{\rm BB}
(\cdot;\hat\sigma,\hat t_1^k,\ldots,\hat t_{\bar\jmath}^k)\bigr)\\
&\quad=\lim_{k\to\infty}\tilde x(1;u^k)
=\lim_{k\to\infty}\tilde\xi^k\\
&\quad=\tilde x^*(1)\\
&\quad=\tilde x\bigl(1;u_{\rm BB}
(\sigma^*,t^*_1,\ldots,t^*_{\tau^*})\bigr),
\end{aligned}
\]
where $\sigma^*:={\rm sgn}(u^*(0))$.
Since the optimal control is (a.e.) unique (Lemma~\ref{LemmaB}), it follows that
\[
\{t_1,\ldots,t_{\bar\imath}\} \supseteq \{t\in(0,1):m_t {\rm~odd}\}
= \{t^*_1,\ldots,t^*_{\tau^*}\}.
\]
Finally, let 
\[
L_i^k:=\min\{\hat t_j^k:i_j=i\}\cup\{t_i\}
\]
where $\cup\{t_i\}$ is not needed when $m_{t_i}$ is even, and
\[
R_i^k:=\max\{\hat t_j^k:i_j=i\}\cup\{t_i\},
\]
then
\[
T_i^k:=[L_i^k,R_i^k] ~~{\rm and~~} |T_i^k|:=R_i^k-L_i^k.
\]
Observe that for all $i$, as $k\to\infty$, $L_i^k\to t_i$, $R_i^k\to t_i$, and
$|T_i^k|\to0$.

\smallskip\noindent
{\em Step 5.}
Remove the finitely many (if any) $\tilde\xi^k$s for which the intervals
$[L_i^k,R_i^k]$, $i=1,\ldots, \bar\imath$ have a nonempty intersection.
Observe that $u^*(t)=u^k(t)$ 
for all $t\in[0,1]\setminus\cup_{i=1}^{\bar\imath}[L_i^k,R_i^k]$ and
let $\Delta u^k:=u^k-u^*$ for all $k$.

Let $I:=\{i: m_{t_i}{\rm ~is~odd}\geq3 {\rm ~and~} t_i\not\in\{0,1\}\}$
and $I^{\rm c}:=\{1,\ldots,\bar\imath\}\setminus I$.
Note that 
$\{t_i\}_{i\in I}\subseteq\{t_i^*\}_{i=1,\ldots,\bar\imath^*}$,
because 
$\{t_i\}_{i\in I}
\subseteq\{t\in(0,1):m_t{\rm ~odd}\}=\{t_1^*,\ldots,t_{\bar\imath^*}^*\}$.
Observe that, for all $i\in I^c$ and all $k\in\N$, $\Delta u^k$ 
(which takes values in $\{0,\pm2\}$) does not 
change nonzero sign on $[L_i^k,R_i^k]$, in the sense that, for all $k$,
$\{-2,2\}\not\subset\{\Delta u^k(t):t\in[L_i^k,R_i^k]\}$.
Indeed, given that the intersection of all intervals $[L_i^k,R_i^k]$ is
empty, if $\Delta u^k$ changes nonzero sign on $[L_i^k,R_i^k]$ for some $k$, 
then $m_{t_i}$ odd, $\geq3$, $t_i\not\in\{0,1\}$, i.e., $i\in I$.
(The fact that $m_{t_i}$ must be odd is a direct consequence of the fact 
that, if $\Delta u^k$ changes nonzero sign on $[L_i^k,R_i^k]$, then $u^*$ 
has to change (nonzero) sign on $[L_i^k,R_i^k]$.)

Next, let $\Delta\tilde x^k(1):=\tilde\xi^k-\tilde x^*(1)$,
$S_i^k:={\rm cl}\{t\in[L_i^k,R_i^k]:u^k(t)-u^*(t)\not=0\}$,
and $|S_i^k|:=\int_{S_i^k}dt$.
Then 
\[
\begin{aligned}
\Delta\tilde x^k(1)&=\int_0^1(u^k(t)-u^*(t))\,\tilde\varphi(t)dt\\
&=\sum_{i\in I^c}\int_{S_i^k}(u^k(t)-u^*(t))
(\tilde\varphi(t_i)+\tilde h_{(i)}(t))dt\\
&\quad+\sum_{i\in I} \int_{S_i^k}(u^k(t)-u^*(t))\\
&\qquad(\tilde\varphi(t_i)+{\rm J}\tilde\varphi(t_i)(t-t_i)\\
&\qquad\quad+\tilde h_{(i)}(t)(t-t_i))dt
\end{aligned}
\]
with $\|\tilde h_{(i)}(t)\|\to0$ as $t\to t_i$ (Taylor, Peano).  Hence
\[
\begin{aligned}
\Delta\tilde x^k(1) &= \sum_{i\in I^c} (w_{i,0}^k\tilde\varphi(t_i) + \tilde r_{(i)}^k)\\
&\quad+\sum_{i\in I}
(w_{i,0}^k\tilde\varphi(t_i)\\
&\qquad+ w_{i,1}^k{\rm J}\tilde\varphi(t_i)\\
&\qquad+ \tilde r_{(i)}^k),
\end{aligned}
\]
where $w_{i,0}^k:=\int_{S_i^k}(u^k(t)-u^*(t))dt$,~
$w_{i,1}^k:=\int_{S_i^k}(u^k(t)-u^*(t))(t-t_i)dt$, and
\[
\tilde r_{(i)}^k:=
\begin{cases}
\int_{S_i^k}(u^k(t)-u^*(t))\tilde h_{(i)}(t)dt & i\in I^c \\
\int_{S_i^k}(u^k(t)-u^*(t))(t-t_i)\tilde h_{(i)}(t)dt & i\in I
\end{cases}
.
\]
Within $S_i^k$, $u^k(t)-u^*(t)$ does not change nonzero sign when 
$i\in I^c$ and  
$(u^k(t)-u^*(t))(t-t_i)$ does not change nonzero sign when $i \in I$. 
The following ensues: 

\noindent
(1) For $i\in I^c$, $|w_{i,0}^k|=\int_{S_i^k} |u^k(t)-u^*(t)|dt$ 
 and for $i\in I$ $|w_{i,1}^k|= \int_{S_i^k} |u^k(t)-u^*(t)| |t-t_i|dt$, and

\noindent
(2) For $i\in I^c$, 
$\|\tilde r_{(i)}^k\|\leq \int_{S_i^k}|u^k(t)-u^*(t)|\|\tilde h_{(i)}(t)\|dt
\leq|w_{i,0}^k| \sup_{t\in S_i^k}\|\tilde h_{(i)}(t)\|$
and for $i\in I$,
$\|\tilde r_{(i)}^k\|\leq \int_{S_i^k}|u^k(t)-u^*(t)|(t-t_i)\|\tilde h_{(i)}(t)\|dt
\leq |w_{i,1}^k| \sup_{t\in S_i^k}\|\tilde h_{(i)}(t)\|$;
hence, for all $i$, $|\tilde r_{(i)}^k\|\to0$ as $k\to\infty$.


\smallskip\noindent
Case $i\in I^c$: Let 
\[
\tilde{\bar h}_{(i)}^k := 
\begin{cases}
\frac{1}{w_{i,0}^k} \tilde r_{(i)}^k & w_{i,0}^k\not=0 \\
0 & {\rm else.}
\end{cases}
\]
Then $\tilde r_{(i)}^k=w_{i,0}^k \tilde{\bar h}_{(i)}^k$ and
$\tilde{\bar h}_{(i)}^k\to0_{1+n}$ as $k\to\infty$.

\smallskip\noindent
Case $i\in I$: Let
\[
\tilde{\bar h}_{(i)}^k := 
\begin{cases}
\frac{1}{w_{i,1}^k} \tilde r_{(i)}^k & w_{i,1}^k\not=0 \\
0 & {\rm else.}
\end{cases}
\]
Then $\tilde r_{(i)}^k=w_{i,1}^k \tilde{\bar h}_{(i)}^k$ and 
$\tilde{\bar h}_{(i)}^k\to0_{1+n}$ as $k\to\infty$.

Now, let $\Delta\tilde x^k(1):= \tilde\xi^k-\tilde x^*(1)$, i.e.,
$\Delta\tilde x^k(1) = \epsilon^k\tilde d+\tilde o^k$,
and let $\tilde h^k:=\frac{1}{\epsilon^k}\tilde o^k$.  Then
\begin{equation}
\label{d+h}
\tilde d + \tilde h^k = (\tilde\Phi+\tilde H^k) \bar w^k,
~~\tilde H^k\to0_{(1+n)\times\bar\ell} {\rm ~as~} k\to\infty,
\end{equation}
where
{\small
\[
\begin{aligned}
\tilde\Phi
&:=\Bigl[(\tilde\varphi(t_i))_{i\in I^c}\quad
({\rm J}\tilde\varphi(t_i))_{i\in I}\\
&\hspace{2.4em}(\tilde\varphi(t_i))_{i\in I}\Bigr]\\
&\in\R^{(1+n)\times\bar\ell},\qquad \bar\ell\leq n
\end{aligned}
\]
and
\[
\begin{aligned}
\tilde H^k:=
[(\tilde{\bar h}_{(i)}^k)_{i\in I^c} ~~ (\tilde{\bar h}_{(i)}^k)_{i\in I}
~~ 0_{(1+n)\times\#I}],\\
\bar w^k:=\frac{1}{\epsilon^k}
[(w_{i,0}^k)_{i\in I^c} ~~ (w_{i,1}^k)_{i\in I} ~~ (w_{i,0}^k)_{i\in I}]^{\rm T}.
\end{aligned}
\]
}
Note that $I$ and $I^c$ satisfy the condition (see theorem statement)
that $\#I^c+3\#I\leq n$.  Indeed
\begin{itemize}
\item $\bar\jmath=\sum_{i=1}^{\bar\imath}m_{t_i}$ by definition of $m_t$;
\item $m_{t_i}\geq3~~\forall i\in I$ by definition of $I$, 
and $m_{t_i}\geq1$ by definition of $I^{\rm c}$;
\item $\bar\jmath\leq n$ by Lemma~\ref{LemmaB},
\end{itemize}
yielding
\[
\begin{aligned}
\# I^c + 3\#I
&= \sum_{i\in I^c}1 + \sum_{i\in I}3\\
&\leq \sum_{i\in I^c}m_{t_i} + \sum_{i\in I}m_{t_i}\\
&= \sum_{i=1}^{\bar\imath}m_{t_i}=\bar\jmath \leq n.
\end{aligned}
\]

\noindent
Equation~(\ref{d+h}) is equivalently written as
\begin{equation}
\label{eq:5g0}
-1+h_0^k = (\Phi_0+H_0^k)\bar w^k
\end{equation}
\begin{equation}
\label{eq:5g1}
\quad 0_n+h^k = (\Phi + H^k)\bar w^k
\end{equation}
where $\tilde\Phi^k = [\Phi_0;\Phi]$ and $\tilde H^k=[H_0;H]$.
From~(\ref{eq:5g0}), $\bar w^k\not\to0$ as $k\to\infty$.

\smallskip\noindent
{\em Step 6.}
Further restrict $\{\xi^k\}$ to an infinite subsequence 
on which $\|\bar w^k\|$ is bounded away from~0.

\smallskip\noindent
{\em Step 7.}
Further restrict $\{\xi^k\}$ to an infinite subsequence where, for some 
$\ell^\times$, argmax$_\ell|\bar w_\ell^k|$ contains $\ell^\times$ at
every $k$.  Then $|\bar w_{\ell^\times}^k|$ is bounded away from 0, and
$\frac{1}{|\bar w_{\ell^\times}^k|}\bar w^k$ is bounded 
(since $\ell^\times$ is selected as an index of maximal absolute value).

\smallskip\noindent
{\em Step 8.}
Restrict $\{\xi^k\}$ to an infinite subsequence (such subsequence
exists, by Bolzano-Weierstra\ss) where 
$\frac{1}{|\bar w_{\ell^\times}^k|}\bar w^k$ converges, say to 
$\underline{\bar w}$; note that $\underline{\bar w}_{\ell^\times}=\pm1$.
Then, dividing both sides of~(\ref{eq:5g1})
by $|\bar w_{\ell^\times}^k|$, invoking~(\ref{d+h}), and taking
the limit as $k\to\infty$ yields $0_n=\Phi\underline{\bar w}$, in contradiction
with $\Phi$ having full column rank.  The proof of the theorem is complete.
\end{proof}

\bigskip
Two results are obtained in the remainder of this subsection.  First, 
we note that, when all eigenvalues of $A$ are real,
the condition in Theorem~\ref{prop:vertical tangent} always holds,
hence Fact~\ref{thm:p0*} holds with $p_0^*=1$ whenever 
$\xi^{\rm f}\in K$.
Then we show that, when some eigenvalues of $A$ are {\em not} real, at
least in the $2\times 2$ case,
Conditions LTS1 still holds when the imaginary parts of such 
eigenvalues are small enough; we also show that, in the cases when
LTS1 fails to hold, it is indeed typically (but not always!) the case 
that Fact~\ref{thm:p0*} only holds with $p_0^*=0$.

As a preliminary, 
note that, as $A$ is real but possibly defective, with only real eigenvalues,
$\varphi^\gamma(t)$ is an ``exponential polynomial'' of the form
\[
\varphi^\gamma(t)=\sum_{i=1}^p P_i(t)e^{a_i t}b_i,
\]
where $a_1,\ldots,a_p$ are distinct real numbers such that
every (possibly repeated) eigenvalue $\lambda$ of $A$ belongs 
to $\{a_1,\ldots,a_p\}$
and each $P_i$ is a polynomial of degree less than the size of the largest
Jordan block associated to eigenvalue $\lambda_i$ of $A$; 
$\varphi^\gamma$ has no more than finitely many zeros unless 
$\gamma=0_n$.  
Vertical tangents can only exist at points $\tilde x(1;u_\gamma)$
where $u_\gamma$ maximizes $x^\gamma(1,u)$  over $u\in\mathcal U$.
Clearly, unless $\gamma=0_n$, the sole such maximizer is the 
only $u_\gamma\in\mathcal U$ 
such that $u_\gamma(t)=$~sgn$(\varphi^\gamma(t))$ for a.e.~$t\in[0,1]$.

\subsubsection{Case when all eigenvalues of $A$ are real}
\label{AwithRealEvalues}


\begin{lemma}
\label{lem:PA2021,Prop3.2}
Suppose $A\in\R^{n\times n}$ and $b\in\R^n$ are such that all eigenvalues
of $A$ are real and $(A,b)$ is controllable. 
Let $t_1<\cdots<t_m$ be in (0,1) and let $i_1<\cdots<i_{\tilde m}$ be in $\{1,\dots m\}$, with $m + \tilde m = n$.
Then 
{\small
\[
\setlength{\arraycolsep}{1pt}
\Phi_n:=\begin{bmatrix} 
e^{(1-t_1)A}b & \cdots & e^{(1-t_m)A}b & A e^{(1-t_{i_1})A}b & \cdots 
& A e^{(1-t_{i_{\tilde m}})A}b \end{bmatrix}
\]
}
is invertible. As a consequence, every matrix obtained by deleting 
columns in $\Phi_n$ has full column rank.  
\end{lemma}
\begin{proof}
``(C1) $\Rightarrow$ (C2)'' portion of 
``Proof of Proposition 4'' in~\cite{AbsilTits2024}, with 
``has zeros at $t_1,\dots,t_n$'' replaced with ``has $n$ real zeros 
counting multiplicities''.
\end{proof}

\noindent
The following theorem is a direct consequence of 
Theorems~\ref{col:xfin K(T_0)}
and~\ref{prop:vertical tangent} and 
Lemmas~\ref{lem:reduced(A,b)} and~\ref{lem:PA2021,Prop3.2}.
\begin{theorem}
\label{thm:PhiNonsingular}
With the same assumptions as in Lemma~\ref{LemmaB},
whenever $\xi^{\rm f}\in K$, 
either (i) $u^*$ is bang-bang or (ii) Fact~\ref{thm:p0*} applied to 
the Lagrange steering problem with LTI dynamics~(\ref{Ax+Bu}) holds 
with $p_0^*=1$ (or both).
\end{theorem}

\subsubsection{Case of $A$ having some eigenvalues off of the real axis}
\label{SomeComplexEvalues}

When some eigenvalues of $A$ are non-real and $\xi^{\rm f}\in\partial K$,
Fact~\ref{thm:p0*} may or may not hold with $p_0^*=1$.
We study the 
situation when $n=2$, i.e., $A:=[\alpha~~\omega;-\omega~~\alpha]$, 
with $\alpha,\omega\in\R$ and, again, assume that $\xi^{\rm f}\in\partial K$,
so that, given $\gamma\in\R^2\setminus\{0\}$, the optimal control
signal $u_\gamma$ is the (a.e.) unique $u_\gamma\in\mathcal U$ that is a maximizer 
in
\[
\gamma^{\rm T}\xi^{\rm f} = \max_{u\in\mathcal U}\gamma^{\rm T} x(1;u).
\]
Accordingly,
$ u_\gamma(t)={\rm sgn}(\varphi^\gamma(t))$ for all $t\in(0,1]$ such
that $\varphi^\gamma(t)\not=0$. 
We have
\begin{equation}
\label{eqn:varphi(t)}
\begin{aligned}
\varphi(t)&=e^{(1-t)A}b \\
&= e^{(1-t)\alpha}
\begin{bmatrix}
b_1\cos\omega(1-t)+b_2\sin\omega(1-t) \\
-b_1\sin\omega(1-t)+b_2\cos\omega(1-t)
\end{bmatrix}
\end{aligned}
\end{equation}
yielding 
\begin{equation}
\label{eq:varphiGamma(t)}
\begin{aligned}
\varphi^\gamma(t) &= e^{(1-t)\alpha}\Bigl(
(\gamma_1b_1+\gamma_2b_2)\cos\omega(1-t)\\
&\qquad +(\gamma_1b_2-\gamma_2b_1)\sin\omega(1-t)
\Bigr)
\end{aligned}
\end{equation}
\[
\begin{aligned}
&=e^{(1-t)\alpha}\Bigl(
(\gamma_1b_1+\gamma_2b_2)(\cos\omega\cos\omega t+\sin\omega\sin\omega t)\\
&\qquad +
(\gamma_1b_2-\gamma_2b_1)(\sin\omega\cos\omega t-\cos\omega\sin\omega t)
\Bigr)
\end{aligned}
\]
\[
\begin{aligned}
&= e^{(1-t)\alpha} (\sigma\cos\omega t + \rho\sin\omega t)\\
&= c(\omega,b,\gamma) e^{(1-t)\alpha}\cos(\omega t+\Theta(\omega,b,\gamma))
\end{aligned}
\]
where
\[
c(\omega,b,\gamma):=\sqrt{\sigma^2+\rho^2}>0, 
\quad \Theta(\omega,b,\gamma):=\arg(\sigma-{\rm j}\rho) ,
\]
with
\[
\begin{aligned}
\sigma&:=(\gamma_1b_1+\gamma_2b_2)\cos\omega+(\gamma_1b_2-\gamma_2b_1)\sin\omega,\\
\rho&:=(\gamma_1b_1+\gamma_2b_2)\sin\omega+(\gamma_2b_1-\gamma_1b_2)\cos\omega.
\end{aligned}
\]
Hence we have, for $t\in(0,1]$,
\begin{equation}
\label{eqn:ugamma=sgn}
u_\gamma(t) = {\rm sgn}(\varphi^\gamma(t))
= {\rm sgn}(\cos(\omega t+\Theta(\omega,b,\gamma)))
\end{equation}
and
\[
\begin{aligned}
\tilde x(1;u_\gamma)&=\int_0^1 \tilde\varphi(t)
{\rm sgn}(\varphi^\gamma(t))dt,\\
\tilde x(1;u_{\gamma,\tilde\epsilon}) &= \int_0^1 \tilde\varphi(t)
~{\rm sgn}(\varphi^\gamma(t)+\tilde\epsilon^{\rm T}\tilde\varphi(t))dt.
\end{aligned}
\]
Now recall that Condition LTS1 holds 
whenever there is no
downward vertical tangent at $\tilde x^*(1)$ and that, for that to
be the case, it is sufficient that, for all $\gamma\in\R^2\setminus\{0\}$,
$\{\varphi(\hat t^j_\gamma), j=1,\ldots,k_\gamma\}$ be a linearly 
independent set, where
$\hat t_\gamma^1<\cdots<\hat t_\gamma^{k_\gamma}$ are the
zeros of $\varphi^\gamma$ in $[0,1]$
(Theorem~\ref{prop:vertical tangent});
also let
$k_\gamma:=0$ when ${\rm sgn}(\varphi^\gamma)$
is constant over $[0,1]$.
When $\omega\in(0,\pi)$, then clearly $k_\gamma=0 {\rm~or~}1$ and it follows
from~(\ref{eqn:varphi(t)}) that no vertical tangent can exist; indeed, 
under our controllability assumption (Lemma~\ref{lem:reduced(A,b)}), i.e.,
$[b_1;b_2]\not=0_2$, $\varphi(t)$ 
cannot vanish anywhere on $\R$ (since $\exp(M)$ is nonsingular for every
matrix $M$); hence, when $\omega\in(0,\pi)$,
Condition LTS1 holds, so that Fact~\ref{thm:p0*} holds 
with $p_0^*=1$, regardless of the location of $\xi^{\rm f}$ in $K$.

{\em Observation.}  Preliminary examination suggests that, in that
$2\times 2$ case, for almost every triple ($\alpha$, $b_1$, $b_2$)
in $\R^3$, when $\omega>2\pi$, Fact~\ref{thm:p0*} holds only with $p_0^*=0$.  
(For $\omega\in(\pi,2\pi)$, the conclusion clearly depends on the data.)

\section{Conclusion}

Conditions 
on finite-horizon optimal control problems without state constraints
were introduced, that imply normality of
Pontryagin's Maximum Principle (i.e., $p_0^*=1$ satisfies the conditions).
The main focus was on the Lagrange problem, for which it was first shown that,
under a regularity condition, whenever the ``lower boundary'' of the
``augmented'' reachable set (i.e., the reachable set of the associated
Mayer problem) is convex and the optimal terminal state lies in the
relative interior of the set of reachable states (of the Lagrange problem), 
the PMP {\em is} normal.  
Next, still within the scope 
of the Lagrange problem, the attention focused on the steering problem
(with target point $\xi^{\rm f}$), with LTI dynamics and a closed
interval of the real line as set $U$ of admissible control values.
There it was shown that (i) when the cost integrand has the form 
$f_0(\xi,\upsilon)=a_0^{\rm T}\xi+h(\upsilon)$, with $h$ continuous,
the PMP is normal 
except, possibly, when the optimal control is bang--bang, 
(ii) when furthermore $h(\upsilon)=b_0\upsilon$, i.e.,
the cost integrand is ``LTI'' as well, {\em and} the eigenvalues of the
state dynamic matrix $A$ are real, 
which is the case, e.g., when $A$ is triangular or symmetric or the system 
is overdamped, then the PMP is {\em always} normal; and, although that is 
still so in some cases where $A$ has some non-real eigenvalues, 
a $2\times 2$ example shows that the condition on the imaginary parts
of the eigenvalues of $A$ cannot be dropped in general.

While the focus of much of the present paper is on single-input problems,
extension to multi-input should be within reach; it is left for future
studies.

\end{document}